\documentclass[11pt]{amsart}
\usepackage{amsfonts, eucal}
\usepackage[all]{xy}
\begin{document}

\newcommand{\C}{{\mathbb C}}                      
\newcommand{\bM}{{\mathcal M}}
\newcommand{\conj}{{\rm conj}}
\newcommand{\mM}{{\mathbb M}}
\newcommand{\Ell}{{\rm Ell}}
\newcommand{\e}{{\rm ell}}
\newcommand{\bEll}{{\bf Ell}}
\newcommand{\half}{{\textstyle{\frac{1}{2}}}}
\newcommand{\cL}{{\mathcal L}}
\newcommand{\dev}{{\rm Dev}}
\newcommand{\oh}{{\mathcal O}}
\newcommand{\pt}{{\rm pt}}
\newcommand{\pairs}{{\rm Pairs}}
\newcommand{\sixth}{{\textstyle{\frac{1}{6}}}}
\newcommand{\T}{{\mathbb T}}
\newcommand{\tate}{{\rm Tate}}
\newcommand{\tmf}{{\rm tmf}}
\newcommand{\Q}{{\mathbb Q}}
\newcommand{\R}{{\mathbb R}}
\newcommand{\UHP}{{\rm \mathfrak h}}
\newcommand{\Sl}{{\rm Sl}}
\newcommand{\PSl}{{\rm PSl}}
\newcommand{\Z}{{\mathbb Z}}
\newcommand{\bLambda}{\mbox{\boldmath$\Lambda$}}

\title{Moonshine elements in elliptic cohomology}

\author{Jack Morava}
\address{Department of Mathematics, The Johns Hopkins University, Baltimore, Maryland
21218}
\email{jack@math.jhu.edu}
\thanks{The author was supported in part by the NSF}
\date {22 November 2007}

\begin{abstract} This is a historical talk about the recent confluence of two
lines of research in equivariant elliptic cohomology, one concerned with {\bf connected} Lie
groups, the other with the {\bf finite} case. These themes come together in (what seems to me
remarkable) work of N. Ganter, relating replicability of McKay-Thompson series to the theory of
exponential cohomology operations. \end{abstract}

\maketitle

\noindent
{\sc Introduction}\bigskip

\noindent
Moonshine, the Monster, conformal field theory, and elliptic cohomology
have all been linked since birth. This survey foregrounds the last of these
topics; its theme is the emergence, under the influence of Hopkins and
Miller and their coworkers, of a  coordinate-free perspective on the moduli 
stacks of bundles over elliptic curves. \bigskip

\noindent
I could never have approached this subject but for the untiring
interest of John McKay, who has often seemed to me an emissary from some
advanced Galactic civilization, sent here to speed up our evolution. I also 
want to acknowledge many helpful conversations about this material with Matt Ando, Andy
Baker, Jorge Devoto, and Nora Ganter. \bigskip

\section{`Classical' elliptic cohomology}

\noindent
{\bf 1.1} The orbifold
\[
\overline{\bM}_\e = (\UHP \cup \Q)/\PSl_2(\Z) 
\]
($\UHP$ is the complex upper half-plane) has one point with isotropy $\Z/2$, one with isotropy
$\Z/3$, and  a third (the cusp) with
(infinite) isotropy $\Z$; it's a good model over $\C$ for the stack $\overline{\bM}_{1,1}(\Z)$ of
stable genus zero curves with one marked point, which is what the real experts work with. These
moduli objects are {\bf not} affine; there is a canonical line bundle $\omega$ over $\overline
{\bM}_\e$, defined by the cotangent line $T^*_{E,0}$ at the origin of the elliptic curve $E \ni
0$. Modular forms are sections of powers of this bundle; they provide a graded
substitute for an affine coordinate ring.\bigskip
 
 \noindent                              
I want to discuss a cohomology theory [23, 30] taking values in the abelian category of
quasicoherent sheaves over $\bM_\e$, such that for example
\[
\bEll(S^{2n}) = \omega^{\otimes n} \;.
\]
Regarding $\omega_E$ as dual to the Lie algebra of $E$ leads to a refinement of this picture
which incorporates odd-dimensional spheres naturally, but that's one of the many topics I'll leave
aside in this talk.\bigskip

\noindent
There is an associated cohomology theory $\Ell^*$ of more classical type, taking values in
modules over a certain graded ring. Away from small primes (two and three) its coefficient ring
is a graded algebra
\[
\Ell^*(\pt) \; \cong \; \oplus_{n \geq 0} \Gamma(\bM_\e, \omega^{\otimes n}) \; \: \: (\otimes
\Z[\sixth])
\]
of classical modular forms, but in general this isomorphism must be replaced by a spectral
sequence
\[
R^*\Gamma(\bM_\e,\omega^{*/2}) \Rightarrow \tmf^*(\pt)
\]
with the derived functor of global sections as its $E_2$ term.  Away from six these
derived functors vanish; but in general the spectral sequence is complicated, with {\bf lots} of
two- and three-torsion, which may look strange to arithmetic geometers but turns out to be  quite
familiar to algebraic topologists. At the prime two, for example, the Hurewicz homomorphism
from the stable homotopy ring of spheres to the resulting cohomology theory is injective up to
dimension fifty or so [25]. \bigskip

\noindent
There are other differences, which leads homotopy theorists to distinguish the abutment of this
spectral sequence, which they call the ring of {\bf topological} modular forms, from its classical
analog. The discriminant $\Delta$, for example, is a modular form but not a topological one: it
supports a nontrivial (torsion-valued) differential in the spectral sequence. However, both
$24\Delta$ and $\Delta^{24}$ survive, and {\bf do} represent elements of $\tmf^*(\pt)$. These
facts seem to be related to some congruences of Borcherds which are beyond  my depth [10, 23
Theorem 5.10];  I mention them only to indicate the arithmetic power of the theory of topological
modular forms.\bigskip

\noindent
There are at present two approaches to constructing topological modular forms, both based a
systematic rigidifications of Landweber's exact functor theorem. The original approach of
Hopkins and Miller builds a homotopy limit of spectra whose connective cover is tmf, while
Lurie's more recent approach constructs this object as the global sections of a sheaf of spectra
over a certain kind of derived moduli stack. The subject is very much under construction (see \S
2.4 below), and both approaches require deep new ideas which change the way we think of both
homotopy theory and algebraic geometry; but in this talk I will take these ideas to some extent
for granted, concentrating instead on their applications to Moonshine and the Monster. \bigskip

\noindent
{\bf 1.2} In practice these constructions have many variants, associated to related moduli problems: 
away from two, for example, curves with $\Gamma_0(2)$ - level structures correspond to Jacobi/Igusa
quadrics 
\[
Y^2 = 1 - 2\delta X^2 + \epsilon X^4 \;,
\]
which was the example which really crystallized this whole subject [22, 29]. Away from
two, the associated modular forms have Fourier coefficients in $\Z[\half]$ and poles only at the
cusp $\tau \rightarrow  i \infty, \; q = e^{2\pi i \tau} \rightarrow 0$, where Tate's elliptic curve
over $\Z[[q]]$ lives. Over the localization $\Z((q))$ its formal group law has multiplicative type,
and restriction to a (formal, {\bf affine}) neighborhood of $i\infty$  defines a natural
transformation 
\[
\bEll \rightarrow K^\tate 
\]
of cohomology theories, which on coefficients becomes the $q$-expansion map
\[
\Ell^* (\pt) \sim \{ {\rm modular \; forms} \} \rightarrow \Z((q)) \;.
\]
The target is a variant of classical $K$-theory, defined by `extension of scalars' from $\Z$ to
$\Z((q))$, suitably oriented: this means that the Chern class of a complex line bundle $L$ in
$K^\tate$ is essentially the Weierstrass $\sigma$-function [4, 34]. \bigskip

\noindent
It is possible to think of $K$-theory as the truncation of $\Ell$ defined by taking the leading term
of a modular form [5, 26, 29]. Some of the more interesting calculations in $K$-theory involve
Bernoulli numbers; in $\Ell$ their analogs involve the Eisenstein series having those numbers as
constant term [7]. The natural operations in $K$-theory are the Adams operations; their analogs
in $\Ell$ are Hecke operators [1,6], as we shall see. \bigskip

\noindent
{\bf 1.3} Physicists (Witten in particular) became interested in these matters because $\Ell^*$ is
the best approximation by a cohomology theory
\[ 
\xymatrix{
{K^*_\T(LM)} \ar[r] \ar[dr] & {\Ell^*(M)} \ar[d]^{q-{\rm exp'n}} \\
& {K^\tate(M)} }
\]
to the $\T$-equivariant $K$-theory of a free loopspace: that construction does not preserve
cofibrations, so the group on the left is very complicated. The completion
\[
u \mapsto 1 - q : K_\T(\pt) = \Z[u,u^{-1}] \rightarrow \Z[[q]] \rightarrow K^\tate(\pt)
\]
restricts to the neighborhood of the fixed point set and so simplifies things greatly [27]; on
the other hand, the resulting theory has a coefficient ring much larger than $\Ell^*$.\bigskip

\noindent
Witten studied the index of an analog of the Dirac operator on $LM$, and showed that when the
loopspace possesses the proper analog of a spin structure, the resulting invariant is a modular
form. This precipitated an enormous amount of research; in particular, such a spin structure is
now understood as the natural form of orientation [4] for elliptic cohomology. Witten further
conjectured the rigidity of his invariant, under actions on the manifold of a {\bf connected} Lie
group $G$. This is now understood in topological terms [3] as a property of the Thom
isomorphism in a suitable form of equivariant elliptic cohomology. \bigskip

\section{The equivariant picture}

\noindent
{\bf 2.1} Sheaf-valued cohomology theories defined over non-affine objects are unfamiliar in
algebraic topology, but working with them is usually routine, perhaps involving higher derived
functors of global sections but not too much more. Extending this framework to the equivariant
context, however, seems to lead us into a genuinely new world: it poses deep {\bf conceptual}
questions. Ganter's recent progress on these questions is the real focus of this talk.\bigskip

\noindent
It seems very reasonable to expect that for a compact Lie group $G$, \medskip

\noindent
$\bullet \;$ there is an equivariant version $X \mapsto \bEll_G(X)$ of elliptic cohomology taking
values in quasicoherent sheaves over some moduli stack 
\[
\bM_\e(G) = \{G-{\rm bundles \; over} \; \check{E} \:|\: E \in  \bM_\e \} \rightarrow \bM_\e
\]
of {\bf principal $G$-bundles} over (the dual of) the `universal' elliptic curve. \bigskip

\noindent
[The appearance here of the dual Abelian variety is technical, having to do with naturality. When
$G = A$ is finite abelian, $A$-bundles over $\check{E}$ are closely related to extensions of
$\check{E}$ by $A$, and thus (by some Cartier duality or Cartan exchange isomorphism) to
homomorphisms from the Pontrjagin dual $\hat{A}$ to $E$; this is related to another approach
[20] to equivariance.]\bigskip

\noindent
{\bf 2.2} For example, when $G = \T$ is the circle group, bundles over $\check{E}$ are
classified by elements of the Picard group
\[
{\rm Pic}(\check{E}) \cong E \;.
\]
A $\T$-space $X$ thus defines a sheaf $\bEll_\T(X)$ over the universal curve ${\bf E} =
\bM_\e(\T)$. \bigskip

\noindent
For {\bf connected} groups this program was pioneered by Grojnowski and 
Ginzburg-Kapranov-Vasserot, followed by Ando and others. For example, the two-sphere $\C
\cup \infty$ with $z \in \T$ acting as multiplication by $z^n$ has (reduced) equivariant elliptic
cohomology defined by the sheaf
\[
\widetilde{\bEll}_\T(S^2(n)) = \oh_{\bf E} (-\sum E[n]) 
\]
[41 \S 3.11]; where $\sum E[n]$ represents the divisor of $n$-torsion points on {\bf E} (equal as
such to $n^2 \cdot [0]$), and $\oh_{\bf E}$ is the structure sheaf of the universal curve. \bigskip

\noindent
The grading of a classical $G$-equivariant cohomology theory $E$ extends naturally from $\Z$
to a Grothendieck group of $E$-orientable linear representations of $G$, but sheaf-valued
cohomology theories over non-affine objects are naturally graded by the Picard group of line
bundles over the base object. When $G$ is connected and (for example) simply-connected,
Grojnowkski [21 \S 3.3]  and Ando [2 \S 9.20, 10.11] use earlier work of Looijenga to construct
a line bundle $\cL$ over $\bM_\e(G)$ such that the $\Z((q))$-module of cusp expansions of
sections of $\cL^{\otimes k}$ can be identified with the corresponding module generated by
graded characters of level $k$ representations of the loop group $LG$.  \bigskip

\noindent
{\bf 2.3} In this talk, however, I want to focus on recent developments in the case of $G$ {\bf
finite}. In spite of deep early work of Devoto, this seems to have received less attention than for
$G$ connected: perhaps because of Witten's rigidity conjecture, and perhaps because of its
intrinsic subtlety. \bigskip

\noindent
We can think of $G$-bundles over $\check{E}$ as classified by elements of the
(zero-dimensional!) orbifold
\[
\pairs_G = {\rm Hom}(\Z^2,G)/G^\conj \cong {\rm Hom}(\pi_1(\check{E}),G)/G^\conj \cong
H^1(\check{E},G)
\]
of conjugacy classes of commuting pairs of elements in $G$ [19]. There is an action of 
$\Sl_2(\Z)$ on such pairs, by 
\[
\left[\begin{array}{cc} a & b \\
                        c & d \end{array}\right] [h,g]
= [h^ag^b, h^cg^d] \;,
\]
and we can take [see [16 \S 2.1] for a more precise statement] 
\[
\bM_\e(G) \sim \pairs_G \times_{\Sl_2(\Z)} \UHP \rightarrow \pt \times_{\PSl_2(\Z)} \UHP \sim
\bM_\e
\]
as a model for the stack of $G$-bundles. Experienced Moonshiners will recognize the space on
the left as a natural habitat for Norton's work. \bigskip 

\noindent
[Note that this map is {\bf not} the identity when $G =\{1\}$ is the trivial
group: the object on the left is then an {\bf un}reduced orbifold, with nontrivial isotropy coming
from the central element $-1$ in $\Sl_2(\Z)$ acting as the involution $[g,h] \mapsto [g^{-1},h^{-
1}]$ on the pair $[1,1]$. This is part of an action of the multiplicative monoid $\Z^\times$ on
pairs by `Adams operations' $[g,h] \mapsto [g^n,h^n]$.] \bigskip

\noindent
Conjugacy classes in the centralizers $C_G(g) = \{ h \in G \;|\: gh = hg \}$ form a  kind of
coordinate atlas on conjugacy classes of pairs, by the map
\[
h \mapsto [h,g] : C_G(g)^\conj \rightarrow \pairs_G
\]
which extends to a parametrization  
\[
C_G(g) \times_\Z ({\rm neighborhood \; of} \; i\infty) \rightarrow \pairs_G \times_{\Sl_2(\Z)}
\UHP = \bM_\e(G)
\]
of the neighborhoods of the cusps in $\bM_\e(G)$; here $\Z$ acts as the subgroup 
\[
\left[\begin{array}{cc} 1 & * \\
                        0 & 1 \end{array}\right]
\]
of translations on pairs, sending $[h,g]$ to $[hg^n,g]$. Restriction to such neighborhoods
of infinity in $\bM_\e(G)$ defines an equivariant lifting
\[
\bEll_G(\pt) \rightarrow K^\dev_G(\pt)
\]
of the $q$-expansion map, $K^\dev$ being Ganter's reformulation [15 \S 3.1] of Devoto's
equivariant version of $K^\tate$, which accomodates some subtle Galois-theoretic properties of
the coefficients of McKay-Thompson series, both classical and generalized [12]. I'll omit some
details here, and describe that construction as the collection, indexed by conjugacy classes of 
elements $g \in G$, of functors of the form 
\[
X \mapsto K^*_{C_G(g)}(X^g)((q^{1/|g|}))|_0 \;,
\]
where $X^g$ is the fixed point set of $g, \; |g|$ is the order of that 
element, and the symbol $|_0$ denotes the
degree zero component of a certain auxiliary $|g|$-grading defined as follows:
\bigskip

\noindent
There are orthogonal projections
\[
 P_k  = |g|^{-1} \sum_{1 \leq n \leq |g|} e^{-2\pi i nk/|g|} g^n \; , \; |g|-1 \geq k \geq 0 \;,
\]
in the complex group ring of the centralizer, which split any $C_G(g)$-equivariant bundle over
$X^g$ into a sum of eigenspaces in which $g$ acts as multiplication by $\exp(2\pi ik/|g)$; this
makes $K^*_{C_G(g)}(X^g)$ into a bigraded algebra (with one grading cyclic of order two, the
other cyclic of order $|g|$). On the other hand $\Z((q^{1/|g|}))$ is also naturally $|g|$-graded, and
the group above carries a tensor product grading. The boundary maps of a cofibration preserve this 
extra structure, so taking its degree zero component is again a cohomology theory. \bigskip

\noindent
At the `classical' cusp (with $g=1$, i.e. associated to pairs of the form $[h,1]$), this construction
simplifies to $K_G(X)((q))$. \bigskip

\noindent
{\bf 2.4} Relations of this form between Ell and $K$-theory are part of a more general picture:
we can think of classical $K_G$ as taking values in sheaves over the moduli space of principal 
$G$-bundles over the circle, which by reasoning like that sketched above is just the space
\[
H^1(S^1,G) \cong {\rm Hom}(\pi_1(S^1),G) \cong G/G^\conj \cong {\rm Spec} \; (R(G) \otimes
\C)
\]
of conjugacy classes in $G$. In fact Quillen showed, back in the 70's, that the classical
cohomology of $G$ can be naturally regarded as a sheaf over the category whose objects are the
abelian subgroups of $G$, with morphisms $A \rightarrow A'$ being the homomorphisms
induced by conjugation by elements of $G$. There is a natural `chromatic' filtration on the
objects of this category (by their rank), with cyclic subgroups at the top. The elliptic picture sees
the first {\bf two} layers of this filtration; the general case is the concern of the generalized
character theory of Hopkins, Kuhn, and Ravenel [24, 35] \bigskip

\noindent
The Hopkins-Miller theorem [19, 39] constructs topological modular forms in terms of a sheaf of
$E_\infty$-algebra spectra over the moduli stack of elliptic curves. Work of Lurie [30 \S 4]
formulates this very naturally in terms of a general representability theorem in enriched (or
derived) algebraic geometry; in particular, he constructs elliptic cohomology as a sheaf of spectra
over a derived moduli stack of suitably oriented elliptic curves over $E_\infty$-algebras. It seems
likely that similar techniques [cf. also [12]] can be used to define an equivariant version of this
theory as a sheaf of spectra over a derived moduli stack, now of $G$-torsors over Lurie's
generalized oriented elliptic curves; but a  precise definition of some such object has yet to
appear.  \bigskip

\section{Ganter's formulation of replicability}
                                                  
\noindent
{\bf 3.1} One concise (but historically misleading) way to tell the Moonshine story is to say that
the graded character of the representation defined by the Frenkel-Lepowsky-Meurman vertex
operator algebra defines a (McKay-Thompson) map $J$ from conjugacy classes in the Monster
to modular functions [9, 11]; on $1 \in \mM$, for example, this construction yields the value
$j(q) - 744$.  In the interpretation presented here, Norton's generalized Moonshine sees $J$ as
the restriction to the classical cusp of a section of $\omega^0$ over $\bM_\e(\mM)$. \bigskip

\noindent
It is known, essentially by case-by-case verification, that the invariance group of a generalized
Moonshine function is a genus zero subgroup of $\Sl_2(\Z)$. This most mysterious property of
McKay-Thompson series is believed (known, in the classical case) to be equivalent to a condition
called {\bf replicability} [11, 36]. The rest of this note is an introduction to recent work of
Ganter, who shows that at the classical cusp, replicability can be expressed very naturally in
terms of cohomology operations on a sheaf-valued theory. \bigskip

\noindent
{\bf 3.2} This involves at least three separate issues, the first being the classification of line
bundles over the $G$-bundle stack. Based on earlier work on Chern-Simons theory, Ganter
[16 \S 2.3] defines a homomorphism 
\[
H^4(BG,\Z) \rightarrow {\rm Pic}(\bM_\e(G)) 
\]
which associates to a degree four cohomology class in $G$, an element of the Picard group of 
line bundles over the moduli space of $G$-bundles. When $G$ is connected and simple, this 
cohomology group is infinite cyclic, and its elements can be naturally identified with the 
{\bf levels} which occur in the theory of positive-energy representations of loop groups [2, 13]. 
On the other hand this cohomology group is finite when $G$ is; for example, there is reason to 
suspect [33 \S 11] that the fourth integral cohomology group of the Monster contains a 
nontrivial element of order 48. \bigskip

\noindent
She then constructs generalizations [16 \S 6.3]
 \[
T_k : \Gamma_\bM(\cL^\alpha) \rightarrow \Gamma_\bM(\cL^{k\alpha}) 
\]
of Hecke operations on the line bundles $\cL^\alpha$ corresponding to $\alpha \in H^4(BG,\Z)$,
which in the case $\alpha = 0$ restrict at the cusp to the classical form
\[
T_k(f(\tau)) = \frac{1}{k} \sum_{ad=k, d > b \geq 0} \psi^a(f)(\frac{a\tau + b}{d}) \;;
\]
where $\psi^a(f)$ denotes the result of applying an Adams operation to the coefficients of the
power series $f$.\bigskip

\noindent
{\bf 3.3} The third topic concerns power operations in generalized cohomology theories. In
$K$-theory the formal sum
\[
V \mapsto \Lambda_t(V) = \sum_{k \geq 0} \Lambda^k(V) t^k : K(X) \rightarrow (1 +
tK(X)[[t]])^\times
\]
of exterior powers of vector bundles defines a homomorphism, because the total exterior power 
satisfies the identity
\[
\Lambda_t(V \oplus W) = \Lambda_t(V) \otimes \Lambda_t(W) \;.
\]
The total symmetric power $S_t(V) = \sum_{k \geq 0} S^k(V)t^k$ is similarly exponential,
in particular because these operations are related by a formal identity
\[
\Lambda_{-t}(V) = S_t(V)^{-1} 
\]
(for example, in some ring of symmetric functions [31]). In that context, Newton's relations lead
to Adams' identity
\[
S_t(V) = \exp(\sum_{k \geq 1} \psi^k(V) \frac{t^k}{k}) \;.
\]
Exponential operations are extremely important in algebraic topology, going back to
work of Atiyah and Steenrod: in modern terms, an orbifold $[X/G]$ has an $n$th orbifold
symmetric power $[X^n/(G \wr \Sigma_n)]$, and a good (multiplicative, equivariant)
cohomology theory $E^*$ will admit `external' power operations
\[
x \mapsto x^{[\otimes n]} : E^*_G(X) \rightarrow E^*_{G \wr \Sigma_n}(X^n) 
\]
which, restricted to the diagonal, yield `internal' power operations: those of Steenrod and 
Adams, as well as Baker's Hecke operations and Ando's higher generalizations (involving
sums over isogenies of formal groups).\bigskip

\noindent
Ganter defined exterior and symmetric power operations for Lubin-Tate theories in her thesis [14
\S 6.11, 7.15], and constructed Hecke operations satisfying the analog [\S 9.2] of the formula
\[
{\bf S}_t(x) = \exp(\sum_{k \geq 1} T_k(x) \frac{t^k}{k}) 
\]
of Newton and Adams. [At about the same time C. Rezk [40 \S 1.12], working in a related but
different context, constructed a kind of universal {\bf logarithmic} operation behaving much like
an inverse to these symmetric powers.] Using the generalized Hecke operations mentioned
above, Ganter extended these constructions to elliptic cohomology and established the \bigskip
  
\noindent
{\bf 3.3 Theorem}[16 \S 6.4]: Replicability of the classical Moonshine function $J$ is equivalent
to the equation
\[
t(J(t) - J(q)) = \bLambda_{-t}(J(q)) \in K^\dev _\mM(\pt)[[t]] \;.
\]

\noindent
{\bf Proof:} The argument uses the theory of Faber polynomials, which associates to
a function
\[
f(q) = q^{-1} + \sum_{k \geq 1} a_k q^k \in \C((q)) \;,
\]
the sequence of polynomials $P_{n,f}(X) = X^n - na_1X^{n-1} + \dots + (-1)^na_n$
characterized by the property
\[
q^{-n} - P_{n,f}(f(q))  \in q\C[[q]] \;.
\]
A generalization [32] of Newton's relations then imply the identity
\[
\log \; \Bigl[q(f(q) - f(p)) \Bigr] = \sum_{n \geq 1} P_{n,f}(f(p)) \frac{q^n}{n} \;;
\]
the assertion of the theorem is thus equivalent to the original form
\[
P_{k,J}(J(q)) = kT_k(J(\tau)) 
\]
of the replicability condition. $\Box$ \bigskip

\noindent
{\bf 3.5} I will close with some remarks, perhaps too vague to be very useful:
\bigskip

\noindent
The first concerns the beautifully simple consequence
\[
q \bLambda_{-t} (J(q)) = - t \bLambda_{-q} (J(t))
\]
of Ganter's formula, reminiscent in many ways of a product formula of Borcherds; 
which, however, at second sight becomes increasingly mysterious. It implies that
the generalized exterior power of $J(q)$ is modular in the purely formal auxiliary
variable $t$. The two sides of the equation are cusp expansions around very different
points: the left-hand side of the equation lies in $K_\mM(\pt)((q))[[t]]$, while the 
right-hand side lies in $K_\mM(\pt)((t))[[q]]$, and their equality implies that both 
lie in $K_\mM(\pt)[[q,t]]$. A coordinate-free interpretation of this, or some similar relation [cf.
[15], just after \S 5.15], might be very enlightening. \bigskip

\noindent                                    
The second is that this replication formula, rewritten as
\[
J(t) = J(q) + t^{-1}\bLambda_{-t}(J(q)) \;,
\]
looks like some kind of evolution equation; but finding an infinitesimal version seems to
be difficult, because of the presence of poles. Related integrable systems will be the subject of 
Devoto's talk.\bigskip

\noindent
The final remark concerns naturality. Once replicability is formulated in terms of cohomology
operations, we can study its behavior under restriction; for example any cyclic subgroup of
$\mM$ pulls $J$ back to a replicable element of $K^\dev_{\Z/n\Z}(\pt)$. [Of course the
resulting element must be the restriction of a Moonshine class for the centralizer of $\Z/n\Z$. The
$n=3$ case, for example, is related to the Thompson group, and is studied in [22]. But even the
pullback to a cyclic group is interesting.]
\bigskip

\noindent
There is evidence [42] for the existence of analogs of $J$ for other sporadic simple groups, and
one might even speculate about symmetric groups. In that context the replicability equation calls
to mind an equation
\[
T(z) = z \cdot \exp (\sum_{k \geq 1} \frac{T(z^k)}{k})
\]
satisfied by the generating function for unlabelled trees (originating with Cayley, but 
in this form due to P\'olya). Current thinking [8] interprets equations of this sort as generalized
{\bf Dyson-Schwinger} equations, taking values in some combinatorially-defined Hopf algebra
(eg, of trees). It is easier to draw a picture to explain where such equations come from, than to
spell it out in words. \bigskip

\noindent
One might speculate that some such equation, in some ring of symmetric functions, might lie
behind the integrable systems investigated in current work of Devoto and McKay. \bigskip

\noindent
Finally, though it seems to have no overt connection to Moonshine, I can't bring myself to end
this talk without mentioning A. Ogg's observation [38] that the primes $p$ dividing the order of
$\mM$ are precisely those for which every supersingular elliptic curve over a finite field of
characteristic $p$ has $j$-invariant in $\mathbb{F}_p$. This seems to be a property neither of
$\mM$ nor of $\bM_\e$ but rather of $\bM_\e(\mM)$; it deserves further attention. \bigskip

\bibliographystyle{amsplain}

\begin{thebibliography}{99}

\bibitem [1]{1} M. Ando, Isogenies of formal group laws and power operations in the
cohomology theories $E_n$, Duke Math. J. 79 (1995) 423--485.

\bibitem[2]{2} ------, Power operations in elliptic cohomology and representations of loop 
groups, Trans. AMS 352 (2000) 5619--5666

\bibitem[3]{3} ------, The sigma orientation for analytic circle-equivariant elliptic cohomology. 
Geom. Topol. 7 (2003) 91--153 

\bibitem[4]{4} ------, M. Hopkins, N. Strickland, Elliptic spectra, the Witten genus and the
theorem of the cube, Invent. Math. 146 (2001) 595--687.
                                                       
\bibitem[5]{5} A. Baker, Operations and cooperations in elliptic cohomology, I: generalized
modular forms and the cooperation algebra, NYJ Math 1 (1994) 39 74

\bibitem[6]{6} --------, Hecke algebras acting on elliptic cohomology, in {\bf Homotopy theory
via algebraic geometry and group representations} 17--26, Contemp. Math AMS (1998)

\bibitem[7]{7} --------, Hecke operations and the Adams $E_2$-term based on elliptic
cohomology, Canad. Math. Bull.  42  (1999) 129--138

\bibitem[8]{8} J.P. Bell, S.N. Burris, K.A. Yeats, Counting rooted trees: the universal law $t(n)
\sim C\\rho^{-n}n^{-3/2}$, Electronic J. Comb 13, Research paper R63 (2006); see 
{\tt  www.combinatorics.org}\

\bibitem[9]{9} R. Borcherds, Vertex algebras, Kac-Moody algebras, and the Monster, Proc. Natl.
Acad. Sci. USA. 83 (1986) 3068--3071

\bibitem[10]{10} --------, Monstrous moonshine and monstrous Lie superalgebras,  Invent. Math. 
109  (1992),  405--444. 

\bibitem[11]{11} C.J. Cummins, T. Gannon, Modular equations and the genus zero property of
moonshine functions, Invent. Math. 129 (1997) 413--443. 

\bibitem[12]{12} J. Devoto, Equivariant elliptic homology and finite groups, Michigan Math. J.
43 (1996) 3--32. 

\bibitem[13]{13} D. Freed, M. Hopkins, C. Teleman, Twisted K-theory and loop group
representations, available at {\tt arXiv:math/0312155}

 \bibitem[14]{14} N. Ganter, Orbifold genera, product formulas and power operations, 
 Adv. Math.  205  (2006) 84--133, available at {\tt  arXiv:math/0407021}
 
\bibitem [15]{15} --------, Stringy power operations in Tate K-theory, available at {\tt 
arXiv:math/0701565} 

\bibitem[16]{16} --------,  Hecke operators in equivariant elliptic cohomology and generalized
moonshine, available at {\tt arXiv:0706.2898}  
                                                                 
\bibitem[17]{17}  D. Gepner, Homotopy topoi and equivariant elliptic cohomology, UIUC
dissertation (2005)

\bibitem[18]{18} V. Ginzburg, M. Kapranov, E. Vasserot, Elliptic algebras and equivariant
elliptic cohomology, available at {\tt arXiv:q-alg/9505012} 
    
\bibitem[19]{19} P. Goerss, M. Hopkins, Moduli problems for structured ring spectra, book in
progress

\bibitem[20]{20} J. Greenlees,  Equivariant formal group laws and complex oriented
cohomology theories, in  {\bf Equivariant stable homotopy theory and related areas},  Homology
Homotopy Appl. 3 (2001) 225--263 

\bibitem[21]{21} I. Grojnowski, Delocalized equivariant elliptic cohomology, in{\bf Elliptic
cohomology, geometry, and applications, and higher chromatic analogs}, LMS Lecture Notes
342, ed. H.R. Miller, D.C. Ravenel, CUP (2007)

\bibitem[22]{22} F. Hirzebruch {\it et al}, {\bf Manifolds and modular forms}, Aspects of
Mathematics, Vieweg (1992)

\bibitem[23]{23} M. Hopkins, Algebraic topology and modular forms, Proceedings of the 2002
ICM, Plenary Lectures and Ceremonies, 291-317 Beijing, Higher Education Press; available
at{\tt arXiv:math/0212397} 

\bibitem [24]{24} ------, N. Kuhn, D. Ravenel, Generalized group characters and complex
oriented cohomology theories,  Jour. AMS 13 (2000)  553--594

\bibitem[25]{25} ------, M. Mahowald: From elliptic curves to homotopy theory, available at {\tt
http://hopf.math.purdue.edu/}

\bibitem[26]{26} N.M. Katz, Higher congruences between modular forms, Ann. Math. 101
(1975) 332--367

\bibitem[27]{27} N. Kitchloo, J. Morava, Thom prospectra for loopgroup representations, in {\bf
Elliptic cohomology, geometry, and applications, and higher chromatic analogs}, LMS Lecture
Notes 342, ed. H.R. Miller, D.C. Ravenel, CUP (2007); available at{\tt arXiv:math/040454}
                                 
\bibitem[28]{28} P.  Landweber, {\bf Elliptic cohomology and modular forms}, Springer LNM
1326 (1988)

\bibitem[29]{29} G. Laures, The topological $q$-expansion principle, Topology 38 (1999) 387--
425
 
\bibitem[30]{30} J. Lurie, A survey of elliptic cohomology, available at {\tt www-
math.mit.edu/}$\sim${\tt lurie/papers/survey.pdf}

\bibitem[31]{31} I. MacDonald, {\bf Symmetric functions and Hall polynomials}, Oxford (1996)

\bibitem[32]{32} J. McKay, The essentials of Monstrous Moonshine, in {\bf Groups and
combinatorics (in memory of Michio Suzuki)} 347--353, Adv. Stud. Pure Math.32 (2001) 

\bibitem[33]{33} G. Mason, Orbifold conformal field theory and cohomology of the Monster, 
available at {\tt http://www.newton.cam.ac.uk/programmes/NST/Mason.pdf}

\bibitem[34]{34} H. Miller, The elliptic character and the Witten genus, in {\bf Algebraic
topology} 281--289 Contemp. Math., 96 AMS (1989)

\bibitem[35]{35} J. Morava,  HKR characters and higher twisted sectors , in {\bf Gromov-Witten
theory of spin curves and orbifolds}  143--152, Contemp. Math., 403, Amer. Math. Soc.,
Providence, RI, 2006,  available at {\tt arXiv:math/0208235}
    
\bibitem[36]{36} S. Norton, More on moonshine, in {\bf Computational group theory (Durham,
1982)} 185--193 Academic Press (1984)
                                                  
\bibitem[37]{37} --------, Appendix to G. Mason, Finite groups and modular functions, in Proc.
Sympos. Pure Math., 47, Part 1, {\bf The Arcata Conference on Representations of Finite
Groups} 181--210 AMS (1987)

\bibitem[38]{38} A. Ogg, Modular functions, in {\bf The Santa Cruz Conference on Finite
Groups} 521--532, Proc. Sympos. Pure Math. 37, AMS (1980)

\bibitem[39]{39} C. Rezk,   Notes on the Hopkins-Miller theorem, in {\bf Homotopy theory via
algebraic geometry and group representations} 313--366, Contemp. Math. 220, AMS (1998)

\bibitem[40]{40} ------, The units of a ring spectrum and a logarithmic cohomology operation,
Jour. AMS 19 (2006) 969--1014, available at {\tt arXiv:math/0407022}

\bibitem[41]{41} I. Rosu, Equivariant elliptic cohomology and rigidity, available at {\tt
arXiv:math/9912089} 
    
\bibitem[42]{42} C. Thomas, {\bf Elliptic cohomology}, Kluwer Academic/Plenum (1999) 

\end{thebibliography}

\end{document}